\documentclass{rspublic} 

\usepackage{graphicx}

\begin{document}

\title{High orders of Weyl series for the heat content}

\author{{Igor Trav\v{e}nec} and {Ladislav \v{S}amaj}}
\affiliation{Institute of Physics, Slovak Academy of Sciences, \\
D\'ubravsk\'a cesta 9, 845 11 Bratislava, Slovakia}


\begin{abstract}{Heat kernel, Heat content, Weyl series}

This article concerns the Weyl series of spectral functions associated 
with the Dirichlet Laplacian in a $d$-dimensional domain with a smooth 
boundary. 
In the case of the heat kernel, Berry and Howls predicted the asymptotic 
form of the Weyl series characterized by a set of parameters.
Here, we concentrate on another spectral function, the (normalized) heat 
content.
We show on several exactly solvable examples that, for even $d$, 
the same asymptotic formula is valid with different values of the parameters. 
The considered domains are $d$-dimensional balls and two limiting cases of 
the elliptic domain with eccentricity $\varepsilon$: 
A slightly deformed disk ($\varepsilon\to 0$) and an extremely prolonged 
ellipse ($\varepsilon\to 1$). 
These cases include 2D domains with circular symmetry and those with only 
one shortest periodic orbit for the classical billiard.
We analyse also the heat content for the balls in odd dimensions $d$ 
for which the asymptotic form of the Weyl series changes significantly.
\end{abstract}

\label{firstpage}
\maketitle

\section{Introduction}
We consider two types of spectral functions: The heat kernel and the 
(normalized) heat content, in the physical literature also known as 
the survival probability. 
They are associated with the spectrum of the Laplacian in a bounded domain. 
Let $\Omega$ be the domain in the $d$-dimensional flat space, with a smooth 
boundary $\partial\Omega$. 
The spectrum of the Laplacian, say with the Dirichlet boundary condition (BC),
is given by
\begin{equation} \label{1}
\begin{array}{rcll}
-\Delta \phi({\bf r}) & = & \lambda \phi({\bf r}) & \qquad{\bf r}\in\Omega , 
\cr 
\phi({\bf r}) & = & 0 & \qquad {\bf r}\in \partial\Omega .
\end{array}
\end{equation} 
The eigenvalues 
$0<\lambda_1<\lambda_2\le \lambda_3 \cdots \le \lambda_j \le \cdots$
form a discrete set (Courant \& Hilbert 1953).
The associated eigenfunctions $\phi_1,\phi_2,\phi_3,\ldots,\phi_j,\ldots$
form an orthonormalized basis of real functions,
$\int_{\Omega} \rd {\bf r} \phi_i({\bf r}) \phi_j({\bf r}) = \delta_{ij}$.

It is instructive to formulate the spectral problem of the Dirichlet Laplacian 
in the context of the diffusion (probability) theory. 
The conditional probability $\rho({\bf r},t\vert {\bf r}_0,0)$ of finding 
a particle at a point ${\bf r} \in \Omega$ at time $t>0$, if it started 
from ${\bf r}_0\in \Omega$ at $t_0=0$, is governed by the diffusion (heat) 
equation
\begin{equation} \label{2}
\frac{\partial \rho({\bf r},t\vert {\bf r}_0,0)}{\partial t} =
\Delta \rho({\bf r},t\vert {\bf r}_0,0) .
\end{equation}
This equation has to be  supplemented by the initial condition
$\rho({\bf r},t=0\vert {\bf r}_0,0) = \delta({\bf r}-{\bf r}_0)$
and by the BC $\rho({\bf r},t\vert {\bf r}_0,0)=0$ for 
${\bf r}\in\partial\Omega$, which reflects absorption of the particle 
hitting the boundary.
The conditional probability can be expressed in terms of the eigenvalues
and eigenfunctions of the Dirichlet Laplacian as follows
\begin{equation} \label{3}
\rho({\bf r},t\vert {\bf r}_0,0) = \sum_j \phi_j({\bf r}_0) \phi_j({\bf r})
{\rm e}^{-\lambda_j t} .
\end{equation}
Various quantities can be constructed from the conditional probability.

Historically, the most studied quantity was the heat kernel, 
defined as the trace
\begin{equation} \label{4}
K(t) = \int_{\Omega} \rd{\bf r}\ \rd{\bf r}_0\ \rho({\bf r},t\vert {\bf r}_0,0) 
\delta({\bf r}-{\bf r}_0)
= \int_{\Omega} \rd{\bf r}\ \rho({\bf r},t\vert {\bf r},0)
= \sum_j {\rm e}^{-\lambda_j t} .  
\end{equation}

Another commonly studied quantity is the heat content, in the chemical 
physics called also the survival probability. 
The ``local'' heat content is defined as
\begin{equation} \label{5}
H(t;{\bf r}_0)= \int_{\Omega} \rd{\bf r} \rho({\bf r},t\vert {\bf r}_0,0) .
\end{equation}
It represents the probability that the particle, localized at a
point ${\bf r}_0\in\Omega$ at $t_0=0$, remains still diffusing
in the domain $\Omega$ at time $t>0$, unabsorbed by the boundary.
The BC is again the Dirichlet one: $H(t,{\bf r}_0) = 0$ 
for ${\bf r}_0\in\partial\Omega$.
The normalized heat content is defined as the average of the local one
over the whole domain,
\begin{equation} \label{6}
H(t) = \frac{1}{\vert\Omega\vert} \int_{\Omega} \rd{\bf r}_0\, H(t;{\bf r}_0)
= \sum_j \gamma_j^2 e^{-\lambda_jt} , \qquad
\gamma_j = \frac{1}{\sqrt{\vert\Omega\vert}} 
\int_{\Omega}\rd {\bf r}\, \phi_j({\bf r}) .
\end{equation}
It represents the probability of finding the particle in $\Omega$ at time 
$t>0$, if it was distributed uniformly with the density $1/\vert\Omega\vert$
over the whole domain at $t_0=0$. 
The normalization by $1/\vert\Omega\vert$ allows us to study also infinite 
domains.
In the mathematical literature, $H(t)$ is interpreted as the amount
of heat at time $t$ inside the domain $\Omega$ with the boundary 
$\partial\Omega$ held at zero temperature, provided that at initial time $0$
the heat was distributed uniformly over $\Omega$.

The heat kernel and the heat content have many similar properties. 
The most important results concern the small-$t$ expansion of both quantities. 
The first three terms of these series are proportional to the volume 
$\vert\Omega\vert$, surface $\vert\partial\Omega\vert$ and 
the integrated mean curvature, respectively. 
The heat kernel $K(t)$ was intensively studied e. g. by Weyl (1911), 
Pleijel (1954), Kac (1966), Stewartson \& Waechter (1971). 
The heat content was introduced by Birkhoff \& Kotik (1954) 
who analysed its 1D version. 
The first few coefficients of the small-$t$ expansion of $H(t)$ for 
Riemannian domains with both Dirichlet and Neumann BC were derived by 
van den Berg et al. (1993), van den Berg \& Gilkey (1994) and 
DesJardins (1998). 
An important progress has been made by Savo (1998a,b) who derived a recurrence 
scheme for the coefficients of the small-$t$ series.

As concerns the heat kernel, it proves useful to apply the Laplace transform 
with respect to $t$. 
For dimensions $d>1$, the heat kernel has to be regularized by subtracting 
the leading terms of the $t$-expansion that would diverge under 
the Laplace integral. 
Thus we get the regularized resolventa $\tilde K(s)$, see e.g. 
Stewartson \& Waechter (1971), Berry \& Howls (1994). 
For example, in the 2D case we have
\begin{equation} \label{8}
\tilde K(s)=\int_0^{\infty} \rd t\, e^{-s^2t} \left( K(t)-
\frac{\vert\Omega\vert}{4\pi t} \right) .
\end{equation}
The small-$t$ expansion of $K(t)$ corresponds to the large-$s$ expansion 
of $\tilde K(s)$:
\begin{equation} \label{9}
\tilde K(s)\sim \sum_{n=1}^\infty \frac{c_n^{(K)}}{s^n} .
\end{equation}
This is the so-called Weyl series.

It turns out that not only the small-$n$ coefficients $c_n^{(K)}$ contain
the information about the domain characteristics. 
Based on a general Borel-transformed theory (Balian \& Bloch 1972, Voros 1983),
Berry \& Howls (1994) developed a formalism for calculating 
the large-$n$ coefficients of the Weyl series in 2D domains. 
They conjectured that
\begin{equation} \label{10}
c^{(K)}_n \sim \alpha \frac{\Gamma(n-\beta+1)}{l^n}, \quad n\to\infty ,
\end{equation}
where $\alpha$ and $\beta$ are some parameters, $l$ is the shortest 
periodic (stable) geodesics for a classical billiard in $\Omega$. 
Note that the factorial nature of the coefficients makes the
Weyl series divergent.
In the case of the disk domain of radius $R$, for which the shortest 
periodic orbits with length $l=4 R$ form a continuous family, the authors 
predicted and tested numerically the parameter values 
$\alpha=1/(8 \sqrt{2\pi})$ and $\beta=3/2$.
In the case of domains with only one (isolated) shortest periodic orbit,
like the ellipse, the authors conjectured the increase of the parameter 
$\beta$ by $1/2$, i.e. $\beta=2$. 

Later Howls \& Trasler (1999) extended these results to higher dimensions. 
They derived exact $\alpha$, $\beta$ and the generalized interpretation 
of $l$ in the case of $d$-dimensional balls of radius $R$. 
They had to distinguish between dimensions $d$ according to their parity. 
For odd $d$, the contribution of the shortest periodic orbit vanishes and 
the next shortest one, i.e. that with three bounces, becomes leading. 
Thus $l$ becomes the perimeter of an equilateral triangle:
\begin{equation} \label{11}
l = \left\{
\begin{array}{ll}
4 R& \mbox{$d$ even,}\cr
3\sqrt{3} R& \mbox{$d$ odd.}
\end{array} \right.
\end{equation}
The parameter $\beta$ depends on $d$ as follows
\begin{equation} \label{12}
\beta = \left\{
\begin{array}{ll}
(5-d)/2& \mbox{$d$ even,}\cr
7/2-d& \mbox{$d$ odd.}
\end{array} \right.
\end{equation}
Explicit formulas for $\alpha(d)$ were also given by Howls \& Trasler 
(1998,1999). 
Further, they proposed a generalization of the asymptotic formula (\ref{10}) 
which involves all periodic orbits $l_j$ on the domain $\Omega$:
\begin{equation} \label{13}
c^{(K)}_n \sim \sum_j \sum_{k=0}^\infty \alpha_{kj} 
\frac{\Gamma(n-\beta_j+1-k)}{l_j^n}, \quad n\to\infty .
\end{equation}

Howls (2001) studied quantum balls threaded by a single magnetic flux
at their center.
While in 2D the parameters $l$ and $\beta$ are insensitive to the periodic 
orbits arising from the diffractive flux line, for a spherical domain these 
parameters are modified by diffractive orbits, in particular $l=2R$ and
$\beta=3$.  

No regularization is needed in the case of the Laplace transform of 
the heat content (van den Berg \& Gilkey 1994). 
We can introduce the Laplace transforms for both the local heat content
\begin{equation} \label{14}
\tilde{H}(s; {\bf r}_0) = \int_0^{\infty} \rd t\, e^{-s^2t} H(t,{\bf r}_0)
\end{equation}
and the heat content
\begin{equation} \label{15}
\tilde{H}(s)=\int_0^{\infty} \rd t\, e^{-s^2t} H(t)
=\sum_{j=1}^{\infty} \frac{\gamma_j^2}{s^2+\lambda_j} .
\end{equation}
These Laplace transforms are related by
\begin{equation} \label{16}
\tilde{H}(s) = \frac{1}{\vert\Omega\vert} \int_{\Omega} \rd {\bf r}_0\,
\tilde{H}(s;{\bf r}_0) .
\end{equation}

The main goal of this paper is the analysis of the counterpart of the
Weyl series (\ref{9}) for the heat content.
The Laplace transform of the heat content has the large-$s$ expansion 
of the form
\begin{equation} \label{17}
\tilde{H}(s)\sim \sum_{n=2}^\infty \frac{c_n}{s^n} .
\end{equation}
We obtain the coefficients $\{ c_n\}$ for few exactly solvable domains. 
For even dimensions $d$, the coefficients $\{ c_n\}$ fulfill at asymptotically 
large $n$ an analogy of the equation (\ref{10}), 
\begin{equation} \label{18}
c_n \sim \alpha \frac{\Gamma(n-\beta+1)}{l^{n-2}}, \quad n\to\infty ,
\end{equation}
where the power $n-2$ in the denominator ensures that $\alpha$ remains 
dimensionless. 
The parameters $\alpha$, $\beta$ and $l$ differ from those for the heat
kernel.
Our task is to determine these parameters for the studied domains
and to point out their possible step-wise modifications under a symmetry 
change of the domain.
A typical example of the symmetry change is a transition from a disk,
possessing an infinite number of the shortest periodic orbits for
the billiard, and an ellipse, with only one shortest periodic orbit.  
For odd dimensions $d$, the Gamma function in the asymptotic equation 
(\ref{18}) has to be replaced by a bounded oscillating function
(at least for the studied $d$-balls).

In connection with Weyl series (\ref{17}) for the heat content we mention 
the paper of van den Berg (2004) who showed that the relation between
the coefficients $c_n$ and the shortest periodic orbit does not hold
in general.
In particular, two domains with different shortest orbits can have
the same $c_n$ coefficients; the shorter of the two periodic orbits
is determined by the difference of the two (exact) heat contents.
This ambiguity of the Weyl series for the heat content inspires us 
to investigate the last for simple domains like $d$-dimensional balls 
and the ellipse.
To our surprise, while a general theory seems to be more complicated for 
the heat content comparing to the heat kernel, the asymptotic Weyl series are 
explicitly available from the exact results and expansions for 
our simple domains.
Similarities and differences between Weyl series for the heat content
and the heat kernel are pointed out.  

The paper is organized as follows. 
In \S 2 we derive a differential equation for $\tilde H(s;{\bf r})$.
In \S 3 we write its exact solution and $\tilde H(s)$ for $d$-balls;
even dimensions are analysed in \S 3a and odd dimensions in \S 3b.
\S 4 is devoted to the 2D ellipse domain of eccentricity $\varepsilon$
and the small-$s$ expansion of $\tilde H(s)$. In \S5, for an arbitrary
value of $s$, we solve exactly the leading terms for two special cases: 
A slightly deformed disk ($\varepsilon\to 0$, \S 5a, \S 5b) and 
an extremely prolonged ellipse ($\varepsilon\to 1$, \S 5c). 
We extract the large-$n$ coefficients of the Weyl series from all exact 
solutions. 
\S 6 is the Conclusion.

\section{Differential equation}
The Laplace transform of the local heat content (\ref{14}) satisfies 
a differential equation which is derived by the following procedure.
Applying the conjugated Laplacian $\Delta^+$ (acting upon the coordinates 
of ${\bf r}_0$) to the representation (\ref{14}) and using the conjugate of 
the diffusion equation (\ref{2}), we obtain
\begin{equation} \label{19}
\Delta^+ \tilde H(s; {\bf r}_0) = \int_0^{\infty} \rd t\, e^{-s^2t} 
\int_{\Omega} \rd {\bf r}\, \Delta^+\rho({\bf r},t\vert {\bf r}_0,0) 
= \int_{\Omega} \rd{\bf r} \int_0^{\infty} \rd t\, e^{-s^2t} 
\partial_t\rho({\bf r},t\vert {\bf r}_0,0) .
\end{equation}
Integration by parts in $t$ then implies the desired equation
\begin{equation} \label{20}
\Delta\tilde H(s; {\bf r}) - s^2 \tilde H(s; {\bf r}) = - 1 ,
\end{equation}
where we abandon the subscript 0 and replace $\Delta^+$ by $\Delta$.
This differential equation is supplemented by the Dirichlet BC
$\tilde H(s;{\bf r}) = 0$ for ${\bf r}\in\partial\Omega$.
A similar differential equation was derived by van den Berg \& Gilkey (1994).

\section{$d$-balls}
For $d$-balls, the Laplacian in (\ref{20}) can be expressed in
terms of $d$-dimensional spherical coordinates.
Since $\tilde H^{(d)}(s;{\bf r})$ does not depend on angle coordinates, 
we end up with the ordinary differential equation
in $r=\vert {\bf r}\vert$; see also van den Berg \& Gilkey (1994).
The solution reads
\begin{equation} \label{22}
\tilde H^{(d)}(s;r) = \frac{1}{s^2} 
\left[ 1-\frac{r^{1-d/2} I_{d/2-1}(s r)}{R^{1-d/2} I_{d/2-1}(s R)} \right] ,
\end{equation}
where $I_{\nu}(x)$ are modified Bessel functions.
Integrating over the domain according to (\ref{16}), we get
\begin{equation} \label{23}
\tilde H^{(d)}(s) = \frac{1}{s^2} 
\left[ 1- \frac{d I_{d/2}(u)}{u I_{d/2-1}(u)} \right] = 
\frac{1}{s^2} \frac{I_{d/2+1}(u)}{I_{d/2-1}(u)},
\end{equation}
where we introduced the scaled variable $u = R s$. 
This result can also be found in the paper of van den Berg \& Gilkey (1994),
though in a slightly more complicated form. 
Its asymptotic analysis depends on the parity of $d$.

\subsection{Asymptotic Weyl series for the balls in even $d$}
Our aim is to find the asymptotic large-$n$ form of the coefficients $c_n$ 
of the Weyl series (\ref{17}) resulting from the equation (\ref{23}). 
It is useful to rewrite the ratio of Bessel functions with the help of 
a logarithmic derivative:
\begin{equation} \label{24}
\tilde H^{(d)}(s) = \frac{1}{s^2}+\frac{d(d-2)}{2u^2 s^2}-\frac{d}{us^2}
\frac{\rd}{\rd u} \ln\ I_{d/2-1}(u) .
\end{equation}
From Gradshteyn \& Ryzhik (2007) we use the large-$u$ expansion
\begin{equation}\label{25}
I_{\nu}(u)\sim \frac{\re^u}{\sqrt{2\pi u}}
\sum_{j=0}^\infty \frac{(-1)^j}{(2u)^j}\frac{\Gamma(\nu+j+
\frac{1}{2})}{\Gamma(j+1)\Gamma(\nu-j+\frac{1}{2})}
\end{equation}
plus exponentially small terms; for the present $\nu=d/2-1$ being integer,
this is an infinite series. 
After straightforward calculation outlined in Appendix A, we get
for asymptotically large $n$
\begin{equation} \label{26}
c_n^{(d)}=(-1)^{d/2-1}\frac{4 d}{\pi\ (2R)^{n-2}}\left[
\Gamma(n-3) + \frac{(d-1)(d-3)}{2} \Gamma(n-4) + \cdots \right] .
\end{equation}

Comparing with the asymptotic representation (\ref{18}) we see that 
the length $l=2R$ is now one half of the shortest periodic orbit $l^{(K)}=4 R$ 
appearing in the case of the heat kernel. 
The symmetry parameter also changes significantly, $\beta=4$ regardless
of (even) dimension. 
Contrary to the heat kernel, the whole dependence on dimension concentrates 
in the prefactor $\alpha^{(d)}=(-1)^{d/2-1}4 d/\pi$.
The next-to-leading term in (\ref{26}) is analogous to the one in
the heat-kernel series (\ref{13}).

\subsection{Asymptotic Weyl series for the balls in odd $d$}
$\tilde H^{(d)}(s)$ in the formula (\ref{23}) with odd $d$ involves Bessel 
functions $I_\nu(u)$ with half-integer $\nu$.
The asymptotic series (\ref{25}) of such functions terminate, 
\begin{equation} \label{28}
I_{\nu}(u)\sim \frac{\re^u}{\sqrt{2\pi u}}
\sum_{j=0}^{\nu-\frac{1}{2}} \frac{(-1)^j}{(2u)^j}\frac{\Gamma(\nu+j+
\frac{1}{2})}{\Gamma(j+1)\Gamma(\nu-j+\frac{1}{2})} , \qquad
u\to\infty .
\end{equation}

For $d=1$ we have
\begin{equation} \label{29}
\tilde H^{(1)}(s) = \frac{1}{s^2}\left(1-\frac{1}{u}\right)
+{\cal O}(1/u^\infty) .
\end{equation}
Consequently, $c_n=0$ for $n>3$. 
An analogous result holds for the Weyl series of the heat kernel. 

For $d=3$, we also find a finite series
\begin{equation} \label{30}
\tilde H^{(3)}(s) = \frac{1}{s^2} 
\left(1-\frac{3}{u}+\frac{3}{u^2}\right)+{\cal O}(1/u^\infty),
\end{equation}
so $c_n=0$ for $n>4$. 
This differs from the heat kernel which has an infinite number of nonzero 
terms in Weyl series for all $d\ge2$, see e. g. Howls \& Trasler (1999).

In the $d=5$ case, the number of non-zero coefficients is infinite: 
\begin{equation} \label{31}
\tilde H^{(5)}(s)\sim \frac{1}{s^2}
\frac{\left(1-6 u^{-1} + 15 u^{-2} - 15 u^{-3} \right)}{1-u^{-1}} 
= \frac{1}{s^2}\left(1-\frac{5}{u}+
\frac{10}{u^2}-5\sum_{k=3}^\infty\frac{1}{u^k}\right) ,
\end{equation}
so that $c_n=-5/R^{n-2}$ for $n\ge 5$. 
We see that the asymptotic formula (\ref{18}) no longer holds:
There is no $\Gamma$-function (thus no $\beta$) in the expression
and the length $l=R$. 

For odd $d>5$, the asymptotic formula for $c_n$ acquires a new 
$n$-dependent term in the numerator. 
In contradiction to the function $\Gamma(n-\beta+1)$, 
this term exhibits bounded oscillations. 
The convergence of the series $\sum_n c_n/s^n$ will depend on $R$.

For $d=7$, these oscillations become periodic. 
Using the asymptotic expansion (\ref{28}) for $I_{5/2}(u)$ 
in the denominator of the representation (\ref{23}) and 
applying the partial fraction decomposition
\begin{equation} \label{32}
\frac{1}{1-3v+3v^2}=\frac{1}{3(v_1-v_2)}\left(\frac{1}{v-v_1}-\frac{1}{v-v_2}
\right)=\sum_{k=0}^\infty \frac{v^k}{3(v_1-v_2)}\left(
\frac{1}{v_1^{k+1}}-\frac{1}{v_2^{k+1}}\right) ,
\end{equation}
where $v=1/u=1/(R s)$ and $v_1 = 3^{-1/2}\exp(\ri\pi/6)$, 
$v_2 = 3^{-1/2}\exp(-\ri\pi/6)$, we get
\begin{equation} \label{34}
c_n=\frac{7}{R^{n-2}}3^{\frac{n}{2}-2}\left[\cos\left(\frac{n\pi}{6}\right)
-\sqrt{3}\sin\left(\frac{n\pi}{6}\right)\right], \qquad n\ge 5 .
\end{equation}
The characteristic length is now $l=R|v_1|=R/\sqrt{3}$. 
The quasi-periodicity $c_{n+12}=3^6 R^{-12} c_n$ is due to 
the commensurability of the phase $\pm\pi/6$ of the complex conjugate 
roots $v_{1,2}$ with $2\pi$.

For $d=9$, we have the polynomial $1-6v+15v^2-15v^3$ in the denominator. 
It has two complex conjugate roots $v_1$, $v_2$ and one real root $v_3$,
expressible in terms of Cardano formulas. 
Numerically, $|v_1|=|v_2|\sim 0.39346201$, with phases 
$\pm 0.2425136494068\ldots \pi$, and $v_3\sim 0.430628846$.
Since the complex phases are no more commensurate with $2\pi$, 
the oscillations of $c_n$ are not periodic. 
The reciprocal of the polynomial $1-6v+15v^2-15v^3$ has three summands of 
the type $v^k/v_m^{k+1},\ \ m=1,2,3$ as in the equation (\ref{32}).
The coefficients $|c_n|$ with large $n$ are dominated by the complex roots 
with the lowest absolute value $|v_1|=|v_2| < v_3$, thus we get $l=R|v_1|$. 
A formal comparison with equation (\ref{13}) leads to the introduction of
two lengths $l_1=R|v_1|$ and $l_2=R v_3$, but we did not find any 
geometric interpretation of these lengths.

For dimensions $d>9$, we get more poles and more summands of type 
$v^k/v_m^{k+1}$. 
The leading term of $|c_n|$ is still given by the pair of conjugate poles 
with the lowest absolute value $v_1$, $v_2=v_1^*$ and the relation $l=R|v_1|$ 
still holds. 
In figure 1, in the units of $R=1$, the roots of the denominator are plotted 
in the complex $v$-plane for larger dimensions $d$. 
The roots with negative imaginary part, placed symmetrically below 
the real axis, are not shown. 
The root $v_1$ with the lowest absolute value is always the leftmost point
of the set. 

To conclude, the consideration of the large-$s$ asymptotic for balls in
odd $d$ leads to $\tilde{H}(s)$ which are the ratios of two polynomials
and as such have a non-zero radius of convergence in $s$.
On the other hand, the Berry and Howls conjecture (\ref{10}) as well as
our conjecture (\ref{18}) have zero radius of convergence in $s$ and
no longer hold.
It is also questionable whether to expand the ratio of two polynomials
into an infinite Weyl series, maybe the geometric information about
the domain is hidden in the polynomial coefficients themselves.  

\begin{figure}
\begin{center}
\includegraphics[scale=0.4]{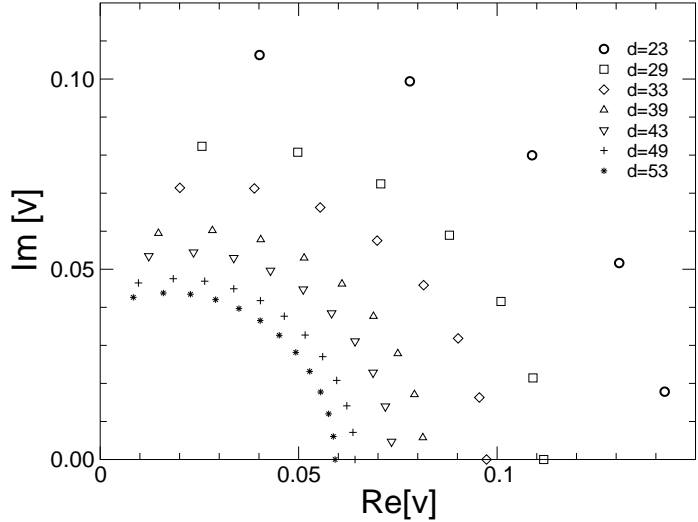}
\caption{Complex poles of $\tilde{H}^{(d)}(s)$ in the $v=1/(Rs)$ plane}
\end{center}
\end{figure}

\section{Ellipse, small-$s$ expansion}
Now, in 2D, we shall pass from the disk to an ellipse which is an example of 
the domain possessing only one shortest periodic orbit. 
It will allow us to reveal a discontinuous change of the $\beta$ 
parameter in (\ref{18}). 

Let us consider the elliptic domain $\Omega$ centered at the origin, 
with major and minor semiaxes $a$ and $b$, respectively. 
In the Cartesian coordinates $(x,y)$, its boundary is given by 
\begin{equation} \label{35}
\frac{x^2}{a^2} + \frac{y^2}{b^2} = 1 , \qquad a\ge b .
\end{equation}
The eccentricity is defined by $\varepsilon = \sqrt{1 - (b/a)^2}\in (0,1)$.
The extreme cases $\varepsilon=0$ and $\varepsilon\to 1$ correspond to
the disk and the infinitely prolonged (locally strip-like) ellipse, 
respectively.

There is little hope to obtain explicitly the heat content for 
the ellipse domain with general $\varepsilon$.
Nevertheless, we are able to construct systematically several types 
of expansions for $\tilde{H}(s;a,b)$. 
We derive the small-$s$ expansion for arbitrary values of $a,\ b$ in 
the present \S4. 
Expansions with respect to $\varepsilon$, around $\varepsilon$ = 0 and 1, 
for arbitrary $s$ will be analysed in the next \S5. 

To find the small-$s$ expansion of $\tilde{H}(s)$ for the elliptic domain,
we first formally expand the local $\tilde H(s;{\bf r})$ in powers of $s^2$,
\begin{equation} \label{36}
\tilde H(s;{\bf r}) = \sum_{j=0}^{\infty} \tilde H_j^{(s)}({\bf r}) s^{2j} .
\end{equation}  
Inserting this expansion into the differential equation (\ref{20}) implies 
an infinite sequence of coupled equations obeyed by the unknown functions 
$\{ \tilde H_j^{(s)}({\bf r})\}$:
\begin{subequations}
\begin{eqnarray}
\Delta \tilde H_0^{(s)}({\bf r}) & = & - 1 , \label{37} \\
\Delta \tilde H_j^{(s)}({\bf r}) & = &  \tilde H_{j-1}^{(s)}({\bf r}) \qquad 
\mbox{for $j\ge 1$.} \label{38}
\end{eqnarray}
\end{subequations}
Each of these functions must satisfy the Dirichlet BC
\begin{equation} \label{39}
\tilde H_j^{(s)}({\bf r}) = 0 \qquad \mbox{for ${\bf r}\in\partial\Omega$,}
\qquad j=0,1,\ldots .
\end{equation}

It is convenient to work with complex coordinates $z=x+iy$ and 
$\bar{z}=x-iy$, in which the elliptic boundary (\ref{35}) becomes
\begin{equation} \label{40}
\frac{1}{2} \left( \frac{1}{a^2} + \frac{1}{b^2} \right) z \bar{z}
+ \frac{1}{4} \left( \frac{1}{a^2} - \frac{1}{b^2} \right) 
\left( z^2 + \bar{z}^2 \right) - 1 = 0 . 
\end{equation}
The Laplacian in complex coordinates has the form
$\Delta = 4 \partial_z\partial_{\bar{z}}$.
We solve successively the set of equations (\ref{37}), (\ref{38})
integrating their r.h.s. in $z$ and $\bar{z}$.
Adding general solutions of the homogeneous equation $\Delta f = 0$, 
\begin{equation} \label{41}
f(z,\bar{z}) = \sum_{j=0}^{\infty} \left( a_j z^j + b_j \bar{z}^j \right) ,
\qquad a_j = b_j ,
\end{equation}
will permit us to fulfill the Dirichlet BC at the boundary (\ref{40}).

Starting with the equation (\ref{37}), we have
\begin{equation} \label{42}
\tilde H_0^{(s)}(z,\bar{z}) = - \frac{1}{4} z \bar{z} + c_0^{(0)} +
c_0^{(1)} (z^2+\bar{z}^2) .
\end{equation}
The coefficients $c_0^{(0)}$ and $c_0^{(1)}$ follow from the condition
$\tilde H_0^{(s)}(z,\bar{z})=0$ at the boundary (\ref{40}):
\begin{equation} \label{43}
c_0^{(0)} = \frac{(ab)^2}{2(a^2+b^2)} , \qquad
c_0^{(1)} = \frac{a^2-b^2}{8(a^2+b^2)} .
\end{equation}
The result (\ref{42}) with (\ref{43}) is substituted into 
the equation (\ref{38}) for $j=1$.
After integrating and adding the homogeneous solution, we obtain
\begin{eqnarray} 
\tilde H_1^{(s)}(z,\bar{z}) & = &
 - \frac{1}{64} (z \bar{z})^2 
+ \frac{1}{4} c_0^{(0)} z \bar{z} + \frac{1}{12} c_0^{(1)} 
( z^3\bar{z} + z\bar{z}^3 ) \nonumber \\ & &
 + c_1^{(0)} + c_1^{(1)} (z^2+\bar{z}^2) + c_1^{(2)}(z^4+\bar{z}^4) . \label{44}
\end{eqnarray}
The coefficients $c_1^{(0)}$, $c_1^{(1)}$ and $c_1^{(2)}$ are again fixed 
to satisfy the Dirichlet BC for $\tilde H_1^{(s)}(z,\bar{z})$ at 
the elliptic boundary. 
We proceed analogously in higher orders. 
The desired expansion of $\tilde H(s)$ in powers of $s^2$ is obtained by 
averaging the relation (\ref{36}) over the ellipse surface:
\begin{equation} \label{45}
\tilde H(s) = \sum_{j=0}^{\infty} \tilde H_j^{(s)} s^{2j} , \qquad
\tilde H_j^{(s)}= \frac{1}{\vert\Omega\vert} 
\int_{\Omega} \rd {\bf r}\, \tilde H_j^{(s)}({\bf r}) ,
\end{equation}
where $\vert\Omega\vert = \pi a b$.
The coefficients of the small-$s$ expansion are obtained in the form
\begin{eqnarray}
\tilde H_0^{(s)} & = & \frac{(ab)^2}{4(a^2+b^2)} \label{46} \\
\tilde H_1^{(s)} & = & - \frac{(ab)^4}{12(a^2+b^2)^2} \label{47} \\
\tilde H_2^{(s)} & = & \frac{(ab)^6[17(a^4+b^4)+98(ab)^2]}{576(a^2+b^2)^3
[(a^4+b^4)+6(ab)^2]} , \label{48} \\
\tilde H_3^{(s)} & = & -\frac{(ab)^8 [93 (a^8+b^8) + 1048a^2b^2 (a^4+b^4)+ 
3190(ab)^4]}{8640(a^2+b^2)^4[(a^4+b^4)+6(ab)^2]^2} , \label{49}
\end{eqnarray}
etc. 

For future purposes, we expand these coefficients around the limit 
$\varepsilon\to 1$, taken as $a \to \infty$ with $b$ fixed. 
We include also the subleading term in $b^2/a^2=1-\varepsilon^2$:
\begin{eqnarray}
\tilde H_0^{(s)} & = & \frac{b^2}{4}-\frac{b^2}{4}\left(\frac{b^2}{a^2}\right)
+{\cal O}\left( \frac{b^4}{a^4}\right) , \nonumber \\
\tilde H_1^{(s)} & = & - \frac{b^4}{12} +\frac{b^4}{6}\left(\frac{b^2}{a^2}\right)
+{\cal O}\left( \frac{b^4}{a^4}\right) , \nonumber \\
\tilde H_2^{(s)} & = & \frac{17}{576}b^6 -\frac{55 b^6}{576}
\left(\frac{b^2}{a^2}\right) +{\cal O}\left( \frac{b^4}{a^4}\right) , 
\nonumber \\
\tilde H_3^{(s)} & = & -\frac{31}{2880} b^8+\frac{11 b^8}{216}
\left(\frac{b^2}{a^2}\right) +{\cal O}\left( \frac{b^4}{a^4}\right) . 
\label{50}  
\end{eqnarray}

\section{Ellipse, eccentricity expansions}
We perform a change of variables $x'= (b/a) x$, $y'=y$;
the Jacobian of this transformation is $J=a/b$.  
The boundary then becomes $x'^2+y'^2=b^2$, i. e. the transformed domain 
$\Omega'$ is the disk of radius $b$. 
The differential equation (\ref{20}) modifies to
\begin{equation} \label{51}
\left[\frac{\partial^2}{\partial {y'}^2} + (1-\varepsilon^2)
\frac{\partial^2}{\partial {x'}^2}\right]\tilde H(s;{\bf r'})
-s^2 \tilde H(s;{\bf r'})=-1.
\end{equation}
Within the probabilistic context explained in Introduction, 
instead of an isotropic diffusion in the ``anisotropic'' ellipse we get 
an anisotropic diffusion in the isotropic disk. 
This approach was inspired by the work Kalinay \& Percus (2006). 

The differential equation (\ref{51}) can be formally expressed as follows
\begin{equation} \label{52}
\left ( \hat{A} + \lambda \hat{B} \right) \tilde{H}(s;{\bf r}')
-s^2 \tilde{H}(s;{\bf r}') = - 1 ,
\end{equation}
where $\lambda$ is a smallness parameter and $\hat{A}, \hat{B}$ are 
the corresponding operators.
There are two natural choices of the smallness parameter $\lambda$. 
In \S5a we shall set $\lambda=\varepsilon^2$ whereas in \S5c 
we shall choose $\lambda=1-\varepsilon^2$. 
In the case $\lambda=\varepsilon^2$, we have
\begin{equation} \label{53}
\hat{A} =\frac{\partial^2}{\partial {y'}^2} +
\frac{\partial^2}{\partial {x'}^2} , \qquad
\hat{B} = -\frac{\partial^2}{\partial {x'}^2} . 
\end{equation}
For $\lambda=1-\varepsilon^2$, we have 
\begin{equation} \label{54}
\hat{A} = \frac{\partial^2}{\partial {y'}^2} , \qquad
\hat{B} = \frac{\partial^2}{\partial {x'}^2} .
\end{equation}
We look for the solution of equation (\ref{52}) perturbatively as 
an infinite series in the smallness parameter $\lambda$:
\begin{equation} \label{55}
\tilde H(s;{\bf r'})= \sum_{n=0}^{\infty} 
\tilde H_n^{(\varepsilon)}(s;{\bf r'}) \lambda^n .
\end{equation}
Inserting this expansion into (\ref{52}) and collecting terms of the same 
powers of $\lambda$, we get a coupled set of differential equations
\begin{eqnarray} 
\hat{A} \tilde{H}_0^{(\varepsilon)}(s;{\bf r'}) - s^2
\tilde{H}_0^{(\varepsilon)}(s;{\bf r'}) & = & -1 , \label{56} \\
\hat{A} \tilde{H}_n^{(\varepsilon)}(s;{\bf r'}) -s^2
\tilde{H}_n^{(\varepsilon)}(s;{\bf r'}) & = & - \hat{B} 
\tilde{H}_{n-1}^{(\varepsilon)}(s;{\bf r'}) \qquad n=1,2,\dots . \label{57}
\end{eqnarray}
All $\tilde{H}_n^{(\varepsilon)}(s;{\bf r}')$ satisfy the Dirichlet boundary 
condition on the disk domain $\Omega'$.
The quantity of interest $\tilde{H}(s)$ is given by
$\tilde{H}(s) = \sum_{n=0}^{\infty} 
\tilde{H}_n^{(\varepsilon)}(s) \lambda^n$, where 
\begin{equation} \label{58}
\tilde H_n^{(\varepsilon)}(s)=\frac{J}{\vert\Omega\vert}\int_{\Omega'}\rd {\bf r'}\, 
\tilde H_n^{(\varepsilon)}(s;{\bf r'}) =\frac{1}{\pi b^2}\int_{\Omega'}
\rd {\bf r'}\, \tilde H_n^{(\varepsilon)}(s;{\bf r'}).
\end{equation}

\subsection{Ellipse, small-$\varepsilon$ expansion}
We first choose $\lambda=\varepsilon^2$ as the smallness parameter. 
Equation (\ref{56}) then becomes
\begin{equation} \label{59}
\left[\frac{\partial^2}{\partial {y'}^2}+\frac{\partial^2}{\partial {x'}^2}
\right]\tilde{H}_0^{(0)}(s;{\bf r'})-s^2 
\tilde{H}_0^{(0)}(s;{\bf r'}) = -1,
\end{equation}
where the upper index $(0)$ refers to the $\varepsilon\to 0$ limit.
The BC is $\tilde{H}_0^{(0)}(s;{\bf r'})=0$ 
for $\vert {\bf r'}\vert = b$.
The solution in polar coordinates $(r',\varphi')$ reads
$\tilde{H}_0^{(0)}(s;r',\varphi')=[1-I_0(sr')/I_0(sb)]/s^2$ which is 
the special case of the equation (\ref{22}) for $d=2$ and $R=b$. 
Applying the relation (\ref{58}) and using the notation $u=s b$, 
we find the 2-ball version of the equation (\ref{23})
\begin{equation} \label{60}
\tilde H_0^{(0)} = \frac{1}{s^2} 
\left[ 1 - \frac{2 I_1(u)}{u I_0(u)} \right].
\end{equation}

The equation (\ref{57}) with $n=1$ takes the form 
\begin{equation} \label{61}
\left[\frac{\partial^2}{\partial {y'}^2}
+\frac{\partial^2}{\partial {x'}^2}\right]
\tilde{H}_1^{(0)}(s;{\bf r'})-s^2 H_1^{(0)}(s;{\bf r'})=
\frac{\partial^2} {\partial {x'}^2}\tilde{H}_0^{(0)}(s;{\bf r'}) .
\end{equation}
In polar coordinates, the r.h.s. becomes
\begin{equation} \label{62}
\frac{\partial^2}{\partial {x'}^2}\tilde H_0^{(0)}(s;r',\varphi')
=-\frac{I_0(sr')}{2I_0(sb)}-\frac{I_2(sr')}{4I_0(sb)}\left(
\re^{2i\varphi'}+\re^{-2i\varphi'}\right) .
\end{equation}
The general solution has the form
\begin{equation} \label{63}
\tilde H_1^{(0)}(s;r',\varphi') = \tilde{h}_0(s;r') +
\tilde{h}_1(s;r') \left( \re^{2i l\varphi'} + \re^{-2i l\varphi'} \right) .
\end{equation}
$\tilde h_0(s;r')$ is determined by the differential equation
\begin{equation} \label{64}
\frac{\rd^2}{\rd {r'}^2}\tilde{h}_0(s;r')+\frac{1}{r'}\frac{\rd}{{\rd} r'}
\tilde{h}_0(s;r') - s^2 \tilde{h}_0(s;r') = - \frac{I_0(sr')}{2I_0(sb)}.
\end{equation}
The homogeneous solutions are $I_0(sr')$ and $K_0(sr')$, 
their Wronskian is $W=-1/r'$. 
Using indefinite integrals like 
$\int z I_0(z)K_0(z) \rd z=z^2[K_1(z)I_1(z)+K_0(z)I_0(z)]/2$ that can be found 
in Bateman \& Erd\'elyi (1953), the solution can be simplified to
\begin{equation} \label{65}
\tilde h_0(s;r')=-\frac{r' I_1(sr')}{4sI_0(sb)}+c_I I_0(sr')+c_K K_0(sr').
\end{equation}
We set $c_K=0$ as we demand regular solution inside the domain, including 
the origin $r'=0$. 
The constant $c_I$ is fixed by the boundary condition $\tilde{h}_0(s;b)=0$.
We thus get
\begin{equation} \label{66}
\tilde h_0(s;r') = - \frac{r' I_1(sr')}{4sI_0(sb)}
+ \frac{b I_1(sb) I_0(sr')}{4sI_0^2(sb)} .
\end{equation}
Since the integral $\int_0^{2\pi} \re^{2i l\varphi'}\rd\varphi'=0$ unless $l=0$, 
only $\tilde h_0(s;r')$ contributes to $\tilde{H}_1^{(0)}$ and 
we finally obtain
\begin{equation} \label{67}
\tilde H_1^{(0)} = - \frac{1}{2 s^2} 
\left[ 1 - \frac{2 I_1(u)}{u I_0(u)} - \frac{I_1^2(u)}{I_0^2(u)} \right].
\end{equation}

Proceeding analogously in higher orders we find
\begin{eqnarray} 
\tilde H_2^{(0)} & = & 
\frac{1}{16 s^2 u^3 I_0^3(u) I_2^2(u)} 
\left[- 2u^3 I_0^5(u)+3u^2 (4+u^2) I_0^4(u) I_1(u) \right. 
\nonumber \\ & & \left. - u(24 + 11u^2)I_0^3(u) I_1^2(u)
+ (16 + 8 u^2 - 3 u^4) I_0^2(u) I_1^3(u) \right. 
\nonumber \\ & & \left.
+ 2 u (2 + 5 u^2) I_0(u) I_1^4(u) - 8 u^2 I_1^5(u) \right] , \label{68} \\
\tilde H_3^{(0)} & = & \frac{1}{192 s^2 u^3 I_0^4(u) I_2^2(u)} 
\left[ u^3 (-12 + 5 u^2) I_0^6(u) \right. 
+ 2 u^2 (36 + 7 u^2) I_0^5(u) I_1(u) \nonumber \\ & & 
- 2u (72 + 57 u^2 + 10 u^4) I_0^4(u) I_1^2(u) 
+ 32 (3 + 4 u^2 + u^4) I_0^3(u) I_1^3(u) \nonumber \\ & & 
+ u (8 + 52 u^2 + 15 u^4) I_0^2(u) I_1^4(u)  
\left. - 36 u^2 (2+u^2) I_0(u) I_1^5(u) + 24 u^3 I_1^6(u) \right] \nonumber \\ 
& & \label{69}
\end{eqnarray}
and so on. 
Recall that $\tilde{H}(s) = 
\sum_{j=0}^{\infty} \tilde{H}_j^{(0)}(s) \varepsilon^{2j}$.

We would like to emphasize that the obtained $\varepsilon$-expansion of
$\tilde{H}(s)$ is valid for {\em all} values of $s$.
This enables us to perform a consistency check of the above results 
by expanding them in small $s$ and comparing with the previous
small-$s$ formulas (\ref{46})-(\ref{49}).
Expanding our $\tilde{H}_j^{(0)}(s)$ in $1/s$ we can also test our results 
by comparison with the exact recurrence scheme of Savo (1998b) for 
the small-$n$ coefficients $c_n$. 
Our results pass also this consistency check.

We can now analyse the large-$s$ behavior of the set (\ref{67})-(\ref{69}), 
in analogy with Appendix A.
Since analytic calculations are cumbersome, they were checked numerically
as well. 
We found that the leading large-$n$ term for the coefficient $c_n$
coming from $\tilde H_k^{(0)}\varepsilon^{2k}$ is proportional to 
$\varepsilon^{2k}[n^k+{\cal O}(n^{k-1})]$. 
Collecting only these leading terms, we get
\begin{equation} \label{70}
c_n(\varepsilon) = \frac{8\Gamma(n-3)}{\pi (2b)^{n-2}}
\left[1-\frac{n\varepsilon^2}{4}+\frac{(n\varepsilon^2)^2}{64}-
\frac{5(n\varepsilon^2)^3}{768}+{\cal O}\left((n\varepsilon^2)^4 
\right)\right] .
\end{equation}
Here, $n\gg N$ where $N$ is a large number. 
So far, the series (\ref{70}) is a formal expansion in $n\varepsilon^2$.
As we are interested in the asymptotically large $n$ for a fixed 
$\varepsilon$, i.e. $n\varepsilon^2\to\infty$, we have to know all terms of 
the series (\ref{70}).
In the next subsection, we propose another method for finding
$\tilde{H}(s)$ based on plausible, but not rigorously justified, arguments.
This method will predict all terms of the expansion (\ref{70}),
reproducing correctly the lowest ones.

\subsection{Renormalized small-$\varepsilon$ expansion}
We return to the original (non-transformed) space ${\bf r}$. 
The ellipse boundary $R(\varphi)$ is expressed in polar coordinates 
as follows
\begin{equation} \label{71}
R(\varphi) = \frac{b}{\sqrt{1-\varepsilon^2\cos^2\varphi}} ,
\qquad \varphi\in (0,2\pi) .
\end{equation}
For small $\varepsilon$, the ellipse is very close to the disk.
Our intuitive approach is based on the assumption that in the
differential equation (\ref{20}) for $\tilde{H}(s;{\bf r})$
we are allowed to neglect the angular part of the Laplacian, i.e.
$\Delta = \partial_r^2 + (1/r)\partial_r$.
The resulting equation is equivalent to that of the disk, the dependence on 
the angle is included only via the elliptic BC at $r=R(\varphi)$:
\begin{equation} \label{72}
\tilde{H}(s;{\bf r}) \sim \frac{1}{s^2} \left[ 1 - 
\frac{I_0(s r)}{I_0(s R(\varphi))} \right] .
\end{equation}
In polar coordinates, the averaging over the ellipse domain is expressible as
\begin{equation} \label{73}
\frac{1}{\vert\Omega\vert} \int_{\Omega} \rd {\bf r} \cdots \equiv
\frac{\sqrt{1-\varepsilon^2}}{\pi b^2} \int_0^{2\pi} \rd\varphi
\int_0^{R(\varphi)} \rd r\, r \cdots .
\end{equation}
After the integration over $r$, we get
\begin{equation} \label{74}
\tilde H(s) \sim \frac{\sqrt{1-\varepsilon^2}}{\pi}
\int_0^{2\pi}\rd\varphi\left[
\frac{1}{2 s^2 (1-\varepsilon^2\cos^2\varphi)}-
\frac{2 I_1(sR(\varphi))}{sb\sqrt{1-\varepsilon^2\cos^2\varphi}
I_0(sR(\varphi))}\right] .
\end{equation}
To analyse the large-$s$ behaviour of the ratio under integration, 
we repeat all steps (\ref{93})-(\ref{96}) of Appendix A, the case $\nu=0$. 
The only difference compared to the 2-ball consists in the replacement of 
the Bessel functions argument $u$ by $sR(\varphi)$:
\begin{equation} \label{75}
\frac{1}{u^j}\to\frac{(1-\varepsilon^2\cos^2\varphi)^{j/2}}{(sb)^j} =
\frac{1}{(sb)^j} \exp\left[{\frac{j}{2}\ln(1-\varepsilon^2\cos^2\varphi)}
\right]\sim\frac{1}{(sb)^j} \re^{-\frac{j\varepsilon^2}{2}\cos^2\varphi} ,
\end{equation}
where we expanded the logarithm only up to the leading $\varepsilon^2$ term. 
The integration over the angle gives
\begin{equation} \label{76}
\int_0^{2\pi} \rd\varphi 
\exp\left(-{\frac{j\varepsilon^2}{2}\cos^2\varphi}\right)
= \exp\left(-{\frac{j\varepsilon^2}{4}}\right)
I_0\left({\frac{j\varepsilon^2}{4}}\right) .
\end{equation}
Similarly as in Appendix A, the $1/s^n$ term in the expansion of
$\tilde{H}(s)$ is identified with the substitution $j=n-4$. 
Using that $e^{-q}I_0(x+q)\sim I_0(x)$ for large $x$, we finally arrive at
\begin{equation} \label{77}
c_n(\varepsilon) \sim\frac{8\Gamma(n-3)}{\pi (2b)^{n-2}} 
\re^{-n\varepsilon^2/4} I_0\left(\frac{n\varepsilon^2}{4}\right) .
\end{equation}

Let us first assume that $n\varepsilon^2$ is finite.
Then
\begin{equation} \label{78}
c_n(\varepsilon) \sim \frac{8\Gamma(n-3)}{\pi (2b)^{n-2}}\sum_{k=0}^\infty
\frac{(2k-1)!!}{(k!)^2}\left(-\frac{n\varepsilon^2}{4}\right)^k .
\end{equation}
The first four terms of this expansion match perfectly those in 
the equation (\ref{70}).
In the limit of interest $n\varepsilon^2\to\infty$, according to (\ref{25}) 
it holds
\begin{equation} \label{79}
c_n(\varepsilon)\sim \frac{8\Gamma(n-3)}{\pi (2b)^{n-2}}
\frac{2}{\sqrt{2\pi n}\ \varepsilon}\sim 
\frac{16\ \Gamma(n-\frac{7}{2})}{\sqrt{2 \pi^3} (2b)^{n-2}\varepsilon} .
\end{equation}
The last expression exploits the property 
$\Gamma(n-7/2)\sim \Gamma(n-3)/\sqrt{n}$ for $n\to\infty$. 
The comparison with the equation (\ref{18}) implies the parameter 
$\beta=9/2$ for the ellipse with small $\varepsilon$. 
This value is by $1/2$ larger than $\beta=4$ of a disk, in close analogy 
with the asymptotic Weyl series for the heat kernel.

For the heat kernel, Berry \& Howls (1994) expected the parameters $\alpha$ 
and $\beta$ to be of the order 1. 
Our result (\ref{79}) suggests something different for the heat content. 
We expect the divergence $\alpha(\varepsilon)\propto 1/\varepsilon$ in 
the symmetry breaking limit $\varepsilon\to 0^+$. 
This is an acceptable price for the step-wise change of $\beta$ 
when restoring the circular symmetry. 
But this is a minor comment, the more important statement about 
the universality of $\beta$ for domains with only one shortest periodic 
orbit still holds.

\subsection{Ellipse with $\varepsilon \to 1$}
Now let us consider $\lambda=1-\varepsilon^2$ as the smallness parameter
in the formalism developed at the beginning of \S 5. 
Equation (\ref{56}) takes the form
\begin{equation} \label{80}
\frac{\partial^2}{\partial {y'}^2}\tilde H_0^{(1)}(s;{\bf r'})
- s^2 H_0^{(1)}(s;{\bf r'})=-1,
\end{equation}
where the upper index $(1)$ refers to the limit $\varepsilon\to 1$. 
From the homogeneous solutions we choose only $\cosh(sy')$, as the odd 
function $\sinh(sy')$ cannot fulfill the boundary condition at 
$y'=\pm\sqrt{b^2-{x'}^2}$. 
Thus we get
\begin{equation} \label{81}
\tilde H_0^{(1)}(s;{\bf r'})=\frac{1}{s^2}\left[1- 
\frac{\cosh(s y')}{\cosh\left(s \sqrt{b^2-{x'}^2}\right)}\right] .
\end{equation}
Now we perform the averaging (\ref{58}), considering for simplicity
four times the first quadrant:
\begin{equation} \label{82}
\tilde H_0^{(1)}(s) = \frac{4}{\pi b^2}
\int_{0}^{b} \rd{x'} \int_{0}^{\sqrt{b^2-{x'}^2}}
\frac{\rd{y'}}{s^2} \left[1- 
\frac{\cosh(s y')}{\cosh(s\sqrt{b^2-{x'}^2})}\right].
\end{equation}
The integration over $y'$ is simple.
To integrate over $x'$, we make the substitution $x'=b \sin\varphi$
and resort to the full angle integration,
\begin{equation} \label{83}
\tilde H_0^{(1)}(s)=\frac{1}{s^2}-\frac{1}{\pi b s^3}\int_0^{2\pi}\cos\varphi
\tanh(b s\cos\varphi)\rd\varphi .
\end{equation}
This is the exact Laplace transform of the heat content for the limiting case 
of an ellipse with finite width $2b$ and infinite length $2a\to \infty$,
valid for any value of $s$. 

To check this expression within the small-$s$ expansion, 
we apply the series representation (Gradshteyn \& Ryzhik 2007)
\begin{equation} \label{84}
\tanh z =\sum_{k=1}^\infty \frac{2^k(2^k-1)}{(2k)!}B_{2k} z^{2k-1},
\end{equation}
where $B_{2k}$ are Bernoulli numbers. 
Inserting this series into (\ref{83}), the integration results in
\begin{equation} \label{85}
\tilde H_0^{(1)}(s)=-2b^2\sum_{k=2}^\infty \frac{(2^{2k}-1)B_{2k}}{(k!)^2} (bs)^{2k-4}
\end{equation}
This series can be compared with the leading terms in the set (\ref{50}) 
and we find a perfect agreement in all available orders. 
It is worth mentioning that for complex $s$ the series (\ref{85}) converges 
if $|s|<\pi/2b=\sqrt{\lambda_1}$, i. e. up to the first
imaginary poles given by the lowest eigenvalue in (\ref{15}).

The next-to-leading term $\tilde{H}_1^{(1)}(s,{\bf r}')$ fulfills the equation
(\ref{57}) with $n=1$:
\begin{equation} \label{86}
\frac{\partial^2}{\partial {y'}^2} \tilde{H}_1^{(1)}(s;{\bf r}')-s^2 
\tilde{H}_1^{(1)}(s;{\bf r}')= -\frac{\partial^2}{\partial {x'}^2}
\tilde{H}_0^{(1)}(s;{\bf r}') .
\end{equation}
With respect to (\ref{81}), the r.h.s. is equal to
\begin{equation} \label{87}
-\frac{\partial^2}{\partial {x'}^2}\tilde{H}_0^{(1)}(s;{\bf r'})
=\cosh(s y')\frac{\partial^2}
{\partial {x'}^2}\frac{1}{s^2\cosh(s\sqrt{b^2-{x'}^2})}
\equiv \cosh(s y') A(x') ,
\end{equation}
where we introduced $A(x')$ for brevity. 
The solution of (\ref{86}) reads
\begin{equation} \label{88}
\tilde{H}_1^{(1)}(s;{\bf r'}) = \frac{A(x')}{2s}y\sinh(sy)-
\frac{A(x')}{2s}\sqrt{b^2-{x'}^2} 
\tanh\left( s\sqrt{b^2-{x'}^2} \right) \cosh(s y') .
\end{equation}
To calculate $\tilde{H}_1^{(1)}(s)$, we first integrate over $y'$ in analogy 
with (\ref{82}), integrate by parts with respect to $x'$ and substitute 
$x'=b \sin\varphi$, to get
\begin{equation} \label{89}
\tilde{H}_1^{(1)}(s) = -\frac{2}{\pi s^2} \int_0^{\frac{\pi}{2}}\rd \varphi 
\sin^2{\varphi}\left[\frac{\sinh^2(bs\cos{\varphi})}{\cosh^4(bs\cos{\varphi})}
+\frac{\sinh^3(bs\cos{\varphi})}
{bs\cos{\varphi}\cosh^3(bs\cos{\varphi})}\right] .
\end{equation}
To check this formula, we expand the integrated function in powers of $s$ 
to get
\begin{equation} \label{90}
\tilde H_1^{(1)}(s)=-\frac{b^2}{4}+\frac{b^4}{6}s^2-\frac{55\ b^6}{576}s^4+
\frac{11\ b^8}{216}s^6-\frac{4487\ b^{10}}{172800}s^8 + \cdots .
\end{equation}
The first four terms can be compared with the $b^2/a^2$ terms 
in the set of four equations (\ref{50}) and we see the full agreement.

Now we are ready to analyse the large-$s$ expansion of the exact solutions 
(\ref{83}) and (\ref{89}). 
The calculations are presented in Appendix B.
Except for the obligatory $1/s^2$ term, only the odd powers of $1/s$ appear 
in the Weyl series.
The results for $c_n$ with $n$ odd are summarized by the formula
\begin{equation} \label{91}
c_n\sim \left[ 1+\frac{\lambda}{2}+
{\cal O}\left(\lambda^2\right) \right] 16\sqrt{\frac{2}{\pi^3}}
\frac{\Gamma(n-7/2)}{(2b)^{n-2}}, \qquad n\to\infty .
\end{equation}
Comparing with the representation (\ref{18}), we see that $l=2b$, 
i.e. one half of the shortest periodic orbit, 
as was expected from the small-$\varepsilon$ analysis. 
The symmetry parameter $\beta=9/2$ is reproduced as well. 
The prefactor is non-universal, dependent on $\lambda= 1-\varepsilon^2$. 

\section{Conclusion}
This paper concerns the asymptotic form of the Weyl series for the
heat content associated with the Dirichlet Laplacian in a smooth
domain $\Omega$.
Using the methods developed by Balian \& Bloch (1972) and Voros (1983),
Berry \& Howls (1994) mapped the quantum billiard model onto the resolventa 
of the heat kernel and conjectured a ``universal'' geometric interpretation 
of the parameters $l$ and $\beta$ in high orders of the Weyl series 
(\ref{10}) {\it for general domains}. 
It is questionable whether an analogical approach to the heat content
is possible.
Some doubts come from the finding of van den Berg (2004) that two
domains with different shortest periodic orbits can have the same Weyl 
series for the heat content. 
One can imagine a scenario analogous to that for the heat kernel where
unstable periodic orbits are excluded from the formalism.
Maybe the accessibility conditions for orbits are even more restrictive for 
the heat content; they might be satisfied for even-dimensional balls and 
the ellipse, but no more for an annulus or twice-cut-disk from van den Berg 
(2004) examples. 
One can also imagine a general analysis of the asymptotic Weyl coefficients
starting from Savo's recurrent scheme (Savo 1998b).

Our strategy was to analyse the asymptotic Weyl series for the heat content,
conjectured in the form (\ref{18}), from the exact results for simple domains.
These results were obtained by solving the differential equation (\ref{20}) 
with Dirichlet BC.
For balls in even dimensions $d$, the conjecture (\ref{18}) applies when
we identify $l=2R$ and $\beta=4$ independently of $d$.
For balls in odd dimensions, $\tilde{H}(s)$ is the ratio of two polynomials
and our conjecture no longer holds.
It might be that the geometric information about the domain is contained
in the polynomial coefficients themselves.
Another open problem is whether the symmetry of balls is more
important than dimensionality, or vice versa, when adapting
our results to non-ball domains. 

Further we studied the ellipse, which represents domains with single
periodic orbit, in two limiting cases of eccentricity $\varepsilon\to 0$ 
and $\varepsilon\to 1$.
In both cases, the parameter $\beta$ is shifted by $1/2$ compared to the disk. 
This phenomenon is in close analogy with the heat kernel.

\begin{acknowledgements}
We are grateful to P. Kalinay for valuable discussions.
This work was supported by the Grants VEGA No. 2/0113/2009 and CE-SAS QUTE.
\end{acknowledgements}

\renewcommand{\theequation}{A.\arabic{equation}}
\setcounter{equation}{0}

\section*{Appendix A}
We are interested in the asymptotic terms of the Weyl series implied by 
the equation (\ref{24}). 
From (\ref{25}), we rewrite the large-$u$ asymptotic of $I_{\nu}(u)$
with $\nu=d/2-1$ as follows
\begin{equation}\label{92}
I_{\nu}(u)\sim \frac{e^u}{\sqrt{2\pi u}} \left( 1 + \Sigma \right) , \qquad
\Sigma = \sum_{j=1}^\infty \frac{(-1)^j}{(2u)^j}\frac{\Gamma(\nu+j+
\frac{1}{2})}{\Gamma(j+1)\Gamma(\nu-j+\frac{1}{2})} .
\end{equation}
Since $\Sigma$ is small in the large-$u$ limit, we can expand 
\begin{equation}\label{93}
\ln I_{\nu}(u) = u -\frac{1}{2}\ln(2\pi u)+\Sigma-\frac{1}{2}\Sigma^2
+ \cdots .
\end{equation}
Further we will explore the identity
$\Gamma(1/2-j) \Gamma(1/2+j) = (-1)^j \pi$ for integer $j$.
In what follows, we shall argue that in the series (\ref{93}) 
$\Sigma$ contributes to the leading and higher-order terms,
$\Sigma^2$ contributes to the subleading and higher-order terms, etc.,
of the coefficients $c_n$ in the asymptotic Weyl series.

Let us first analyse just the leading term given solely by $\Sigma^1$.
Substituting the truncated expansion (\ref{93}) into (\ref{24})
and using standard properties of $\Gamma$-functions, we get 
the leading contribution
\begin{equation}\label{94}
\tilde H^{(d)}(s) \sim \frac{d}{u s^2}\sum_{j=1}^\infty \frac{1}{2^ju^{j+1}} 
\frac{(\frac{d-3}{2}+j)...(\frac{3}{2}+j)(\frac{1}{2}+j)}{(\frac{d-3}{2}-j)
(\frac{d-5}{2}-j)... (\frac{1}{2}-j)}
\frac{\Gamma^2(\frac{1}{2}+j)}{\pi\Gamma(j)} .
\end{equation}
Further we need the large-$j$ expansion
\begin{equation}\label{95}
\frac{\Gamma^2(\frac{1}{2}+j)}{\Gamma(j)\Gamma(j+1)}\sim 1-\frac{1}{4j}
+\cdots \quad j\gg 1 .
\end{equation}
Taking the limit $j\to\infty$ of the long ratio in (\ref{94}) 
we get $(-1)^{d/2-1}$. 
Considering the leading (first) term of the expansion (\ref{95}), we obtain
\begin{equation}\label{96}
\tilde H^{(d)}(s) \sim (-1)^{d/2-1}\frac{d}{R s^3}
\sum_{j=1}^\infty \frac{1}{2^j (Rs)^{j+1}} \frac{\Gamma(j+1)}{\pi} .
\end{equation}
Comparing this series with (\ref{17}), we set $j=n-4$ and find 
the leading term of the asymptotic formula (\ref{26}).

The calculation of the subleading $1/n$ term in the equation (\ref{26}) 
is more complicated and, for simplicity, we restrict ourselves
to the $d=2$ case. 
There are two contributions. 
The simpler one comes from the subleading term in (\ref{95}), i.e. we get 
the leading factor $8\Gamma(n-3)/\pi(2R)^{n-2}$ times $(-1/4n)$. 
The tricky part comes from the $\Sigma^2$ term in (\ref{93})
which contributes to $c_{n+2}/[s^2(2u)^{n}]$ by the sum
\begin{eqnarray} \label{97}
- \sum_{k=1}^{n-1}\frac{\Gamma^2(\frac{1}{2}+k)}{2\pi^2\Gamma(k+1)}
\frac{\Gamma^2(\frac{1}{2}+n-k)}{\Gamma(n-k+1)} & = &
\left[ 2 - 
{_3F_2}\left( \frac{1}{2},\frac{1}{2},-n;\frac{1}{2}-n,\frac{1}{2}-n;-1
\right) \right] \nonumber \\ & & \times
\frac{\Gamma^2(\frac{1}{2}+n)}{2\pi\Gamma(1+n)} .
\end{eqnarray}
Here, we introduced the hypergeometric function $_3F_2$, see e.g. 
Gradshteyn \& Ryzhik (2007). 
Its large $n$ asymptotic is
\begin{equation}\label{98}
_3F_2\left(\frac{1}{2},\frac{1}{2},-n;
\frac{1}{2}-n,\frac{1}{2}-n;-1\right)\sim 2
+\frac{1}{2n}+{\cal O}\left(\frac{1}{n^2}\right), \qquad n\to\infty .
\end{equation}
Inserting this into (\ref{97}), in equation (\ref{24}) 
we get exactly the same contribution as above, i.e. the leading factor times
$(-1/4n)$. 
Summing up two equal contributions we get finally $8\Gamma(n-3)/\pi(2R)^{n-2}$ 
times $(-1/2n)$; noting that $\Gamma(n-4)\sim\Gamma(n-3)/n$ for large $n$,  
this is already the subleading term of the equation (\ref{26}) for 
the special case $d=2$. 
One can show that $\Sigma^3$, etc. contribute to higher-order terms.

The derivation can be generalized to higher dimensions $d$.

\renewcommand{\theequation}{B.\arabic{equation}}
\setcounter{equation}{0}

\section*{Appendix B}
We aim at analysing the asymptotic Weyl expansion of the equation (\ref{83}). 
Let us return to the first quadrant integration and rewrite appropriately 
$\tanh$, 
\begin{equation} \label{99}
\tilde H_0^{(1)}(s)=\frac{1}{s^2}-\frac{4}{\pi b s^3}
+\frac{8}{\pi b s^3}\int_0^{\frac{\pi}{2}}
\frac{\cos\varphi \rd\varphi}{e^{2bs\cos\varphi}+1}
\end{equation}
Further we apply the substitution $z=2bs\cos\varphi$ and subsequently 
the series
\begin{equation} \label{100}
\frac{q}{\sqrt{1-q^2}}=\sum_{k=0}^\infty\frac{(2k-1)!!}{k!2^k}q^{2k+1} .
\end{equation}
We need to calculate the integral $\int_0^{2bs} \rd z\, z^{2k+1}/(e^z+1)$. 
For very large $s$, the upper integration limit can be extended to infinity, 
because the added term is exponentially small
\begin{equation} \label{101}
0\le\int_{2bs}^\infty\frac{z^{2k+1}\rd z}{e^z+1} <\int_{2bs}^\infty\frac{z^{2k+1}
\rd z}{e^z}=(2k+1)!e^{-2bs}\sum_{m=0}^{2k+1}\frac{(2bs)^m}{m!} .
\end{equation}
We can calculate the large-$k$ asymptotic of
\begin{equation} \label{102}
\int_0^\infty\frac{z^{2k-1}\rd z}{e^z+1}\sim (2k-1)! \qquad k\to\infty .
\end{equation}
Further we use the obvious relation $(2k-1)!!=(2k)!2^{-k}/k!$. 
Applying the above steps, we get for both $s$ and $k$ large
\begin{equation} \label{103}
\tilde H_0^{(1)}(s)\sim \frac{8}{\pi}\sum_k \frac{(2k)!(2k+2)!}{
2^{2k}k!(k+1)!}\frac{1}{s^2(2bs)^{2k+3}}, \qquad k\to\infty .
\end{equation}
Considering the asymptotic behaviour of
\begin{equation} \label{104}
\frac{\Gamma(2k+1)\Gamma(2k+3)}{2^{2k}\Gamma(k+1)\Gamma(k+2)}\sim 
2\sqrt{\frac{2}{\pi}}\Gamma\left(2k+\frac{3}{2}\right)
\left[1+{\cal O}\left(\frac{1}{k}\right)\right], \qquad k\to\infty
\end{equation}
and identifying $n=2k+5$, we arrive at the first term 
in the equation (\ref{91}). 
Note that $n$ is odd, the even powers of $1/s$ do not appear.

\subsection*{Subleading term in $\lambda=1-\varepsilon^2$}
Let us first rewrite the equation (\ref{89}) in the following way
\begin{equation}\label{105}
\tilde H_1^{(1)}(s) = -\frac{2}{3\pi s^2}\frac{\rd}{\rd s} F(s) , \qquad 
F(s) = s^4\int_0^\frac{\pi}{2}
\rd \varphi \sin^2{\varphi}\frac{\tanh^3{w}}{w} ,
\end{equation}
where $w=bs\cos{\varphi}$. 
$F(s)$ fulfills the relation
\begin{equation}\label{106}
\frac{\rd}{\rd s}\left[s F(s)\right]=3\int_0^{b s}\frac{\rd w}{b s}
\sqrt{1-\left(\frac{w}{b s}\right)^2}
\frac{\sinh^2{w}}{\cosh^4{w}}.
\end{equation}
We expand
\begin{equation}\label{107}
\sqrt{1-\left(\frac{w}{b s}\right)^2}\frac{\sinh^2{w}}{\cosh^4{w}}\sim
-\sum_{k=0}^\infty \frac{(2k-2)!}{2^{2k-1}k!(k-1)!}
\left(\frac{w}{b s}\right)^{2k}\left[4\re^{-2w} + {\cal O}(\re^{-4w})\right].
\end{equation}
The upper limit of the integration in (\ref{106}) can be again 
extended to infinity, as in the case of the equation (\ref{101}).
Then we calculate $[s F(s)]'$ and after all we get
\begin{equation}\label{108}
F(s)\sim 12\sum_k\frac{(2k-2)!(2k-1)!}{2^{2k-1}k!(k-1)!}\frac{1}{(2bs)^{2k+1}}.
\end{equation}
We substitute this series into (\ref{105}) and use the large-$k$ behavior of 
the ratio
\begin{equation}\label{109}
\frac{(2k-3)\Gamma(2k-1)\Gamma(2k)}{2^{2k-1}\Gamma(k+1) \Gamma(k)}\sim
\sqrt{\frac{2}{\pi}}\Gamma\left(2k-\frac{1}{2}\right) , \qquad k\to\infty.
\end{equation}
Thus we obtain
\begin{equation}\label{110}
\tilde H_1^{(1)}(s)\sim\frac{8}{\pi}\sqrt{\frac{2}{\pi}}\sum_k
\frac{\Gamma\left(2k-\frac{1}{2}\right)} {s^2(2bs)^{2k+1}} , \qquad k\to\infty.
\end{equation}
We set $n=2k+3$ and finally arrive at the second term in 
the equation (\ref{91}). 
Note that only odd powers of $1/s$ appear again.

\end{document}